\documentclass[oneside]{amsart}

\usepackage[english]{babel}
\usepackage{amsmath}
\usepackage{amsthm}

\newtheorem{theorem}{Theorem}[section]
\newtheorem{corollary}[theorem]{Corollary}
\newtheorem{proposition}[theorem]{Proposition}
\newtheorem{lemma}[theorem]{Lemma}

\theoremstyle{definition}
\newtheorem{definition}[theorem]{Definition}
\newtheorem{example}[theorem]{Example}
\newtheorem{remark}[theorem]{Remark}

\usepackage{amssymb}
\usepackage{mathtools}
\usepackage[unicode=true,hidelinks,bookmarksnumbered]{hyperref}
\usepackage{enumitem}

\setlist{noitemsep}
\setlist[enumerate]{label=(\arabic*), font=\upshape}

\usepackage{tikz}

\usetikzlibrary{babel}
\usetikzlibrary{matrix}
\usetikzlibrary{arrows}
\usetikzlibrary{calc}
\usetikzlibrary{cd}
\usetikzlibrary{positioning}
\usetikzlibrary{decorations.markings}
\usetikzlibrary{svg.path}
\tikzset{ampersand replacement=\&}

\usepackage[style=numeric,maxbibnames=99,maxalphanames=5]{biblatex}
\usepackage{csquotes}

\DeclareMathOperator\skel{skel}
\DeclareMathOperator\colim{colim}
\DeclareMathOperator\Ker{Ker}

\DeclareMathOperator\Hom{Hom}

\DeclareMathOperator\LS{LS}
\DeclareMathOperator\ug{u}
\newcommand\lex{\mathrm{lex}}

\newcommand\ord{\mathrm{ord}}
\newcommand\Simp{\operatorname{Simp}}
\newcommand\Set{\operatorname{Set}}
\newcommand\Mod[1]{\operatorname{Mod}{#1}}
\newcommand\pt{\mathrm{pt}}

\begin{document}
\title[Cofiltrations of spanning trees]{Cofiltrations of spanning trees in multiparameter persistent homology}
\author{Fritz Grimpen}
\email{grimpen@uni-bremen.de}
\address{Institute for Algebra, Geometry, Topology and its Applications (ALTA), Department of Mathematics, University of Bremen, Germany}
\author{Anastasios Stefanou}
\email{stefanou@uni-bremen.de}
\address{Institute for Algebra, Geometry, Topology and its Applications (ALTA), Department of Mathematics and Computer Science, University of Bremen, Germany}
\begin{abstract}
    Given a multiparameter filtration of simplicial complexes, we consider the problem of explicitly constructing generators for the multipersistent homology groups with arbitrary PID coefficients.
    We propose the use of spanning trees as a tool to identify such generators by introducing a condition for persistent spanning trees, which is accompanied by an existence result for cofiltrations consisting of spanning trees.
    We also introduce a generalization of spanning trees, called spanning complexes, for dimensions higher than one, and we establish their existence as a first step towards this direction.
\end{abstract}
\maketitle

\tableofcontents

\section{Introduction}
Persistent homology studies the evolution of topological features of a dataset (finite metric space) across a given filtration on the dataset \cite{EdelsbrunnerLetscherZomorodian2002}.
It is arguably the most heavily used tool in Topological Data Analysis (TDA) \cite{Carlsson2009,EdelsbrunnerHarer2010}.
One-parameter persistent homology has the algebraic structure of an $\mathbb Z$-graded module \cite{CarlssonZomorodianCollinsGuibas2005}.

\subsection{Multiparameter persistence} Persistent homology extends naturally from one-parameter filtrations \cite{ZomorodianCarlsson2005} to multiparameter filtrations of datasets \cite{CarlssonZomorodian2009}. The corresponding algebraic structure of multiparameter homology is that of an $\mathbb Z^d$-graded module, as shown in \cite{CarlssonZomorodian2009}.
In fact persistent homology can be defined for filtrations parameterized over any ordered set $Q$ \cite{BubenikScott2014}.
The extension from graded modules to multigraded modules is of particular importance to TDA, since there are situations where important information from data cannot be captured from persistent homology by merely considering one-parameter filtration, e.g.~persistent homology is shown to be unstable to the presence of outliers \cite{CarlssonZomorodian2009,LesnickWright2015}.
By introducing a density parameter as a second filtration parameter on datasets, outliers in datasets can then be regarded as noise, making two-parameter persistent homology robust to the presence of outliers in the data. Thus, additional information can be extracted from datasets by studying multiparameter persistent homology.
In multiparameter persistence, in general, there is no way to realize the indecomposable modules as interval modules \cite{CarlssonZomorodian2009, DeyXin2022}, as in the case of one-parameter persistence \cite{BotnanCrawley-Boevey2020,Crawley-Boevey2015}.
Therefore, recent research works in multipersistence theory have been focusing on the study of multipersistent invariants \cite{KimMemoli2021,Patel2018}.

\subsection{Minimal presentations}
The key techniques to better understand the structure of multiparameter persistence modules as well as to help towards an efficient computation of multipersistent invariants of modules, are the minimal free presentations, minimal free resolutions \cite{CoxLittleOShea2005}, and the recently introduced  minimal flat-injective presentations \cite{Miller2020a} of multigraded modules.
However, in practice, all of these techniques often pre-assume that the module is equipped with a set of generators to begin with, and that the multigraded module is a submodule of some given free module.
Gröbner basis methods can then be utilized to compute a minimal presentation and a minimal free resolution of the module \cite{CoxLittleOShea2005,LaScalaStillman1998,MillerSturmfels2005}.
As of now, the problem of finding an efficient way to compute a minimal presentation of multiparameter persistence modules remains open.
Some works in TDA towards this direction include \cite{CarlssonSinghZomorodian2010,ChachólskiScolamieroVaccarino2017,LesnickWright2022}.

\subsection{Related work}
In \cite{CarlssonSinghZomorodian2010} the authors show that if a multifiltration $X = (X_a)_{a \in \mathbb N^d}$ is one-critical, i.e.~every simplex enters the multifiltration in at most one multidegree (thereby making each $C_n(X)$ a free module), then Schreyer's algorithm \cite{CoxLittleOShea2005} can be directly applied for computing a Gröbner basis for the $n$-cycle module $Z_n(X)=\Ker \partial_n$ of the homology module $H_n(X)=\Ker \partial_n / \operatorname{Im}\partial_{n+1}$ of a multifiltration $X$, and thus for computing a minimal free presentation for the homology module.
In \cite{ChachólskiScolamieroVaccarino2017}, the authors developed an explicit and combinatorial description for the homology modules of a multifiltration of simplicial complexes.
They developed a construction that realizes the $n$-th homology module $H_n(X)$ of a multifiltration $X$, as the homology  $\Ker f/\operatorname{Im} g$ of a certain chain complex  $F_1\xrightarrow{f}F_2\xrightarrow{g}F_3$ of free multigraded modules.
They also proved that the multigraded modules that can occur as $k$-spans of multifiltrations of sets are the direct sums of monomial ideals, where $k$ is an arbitrary commutative, unitary ring.

In the case of a $2$-parameter filtration $X$, they construct a free presentation
\[ F_1\to F_0\to H_n(X)\to 0 \]
for the homology module $H_n(X)$, where $F_1$ is explicitly defined and $F_0$ is the kernel of a  morphism between certain free modules.
In \cite{LesnickWright2022} the authors are divising a general algorithm for computing minimal presentation for $2$-parameter persistence modules (also called bigraded modules), whenever they are represented as $M\cong \Ker g/\operatorname{Im} g$ (for instance via the combinatorial presentation of \cite{ChachólskiScolamieroVaccarino2017}), for some chain complex of free $2$-parameter persistence modules $F_1\xrightarrow{f}F_2\xrightarrow{g}F_3$.
In both works, for $d=2$, the authors are invoking the algebraic-geometric fact that over two variables ($2$-parameters), kernels of free bigraded modules are also free.
This is no longer true when we work on multigraded modules over $d \geq 3$ parameters.

\subsection{Motivation}
Our goal is to develop a method for computing (minimal) generators for the homology of a multifiltration.
Our work is motivated by the problem of computing generators of kernels of morphisms of \emph{upper set decomposable} modules (direct sums of monomial ideals) as in \cite{ChachólskiScolamieroVaccarino2017}.
In this work our approach is to investigate a way to find (minimal) free generators for the $n$-cycle module $Z_n(X)$ of a multifiltration $X$, since this would be enough for our goal.
Indeed, if we know the generators of $Z_n(X)$, because of the canonical epimorphism from $Z_n(X)$ onto $H_n(X)$, we will also obtain free generators for $H_n(X)$.
This will then enable us, at least implicitly, to be able to apply the machinery of commutative algebra to compute a minimal free presentation and a minimal free resolution of the homology of a multifiltration.
In the static case of simplicial complexes, we can compute the generators of the first homology group $H_1(X)$ of a given complex $X$ by using spanning tree subcomplexes $T$ of $X$.
We have the following theorem which is known from many sources in the literature \cite{Kozlov2020,Spanier1982}.

\begin{theorem}
    If $X$ is a 1-dimensional simplicial complex and $T$ a spanning tree of $X$, then $H_1(X) \cong C_1(X, T)$.
\end{theorem}

\subsection{Our contributions}

The main contribution of this paper is a method to construct generators of the $1$-dimensional cycle module $Z_1(X; k) = \Ker \partial_1$ for a given multifiltration $X = (X^q)_{q \in Q}$ of finite simplicial complexes over a finite ordered set $Q$ and an arbitrary principal ideal domain~$k$.
To the multifiltration $X$ we associate a \emph{cofiltration of spanning trees} $T = (T^q)_{q \in Q}$ by endowing the simplicial complex $\bigcup_{q \in Q} X^q$ with a simplicial order.
We show that these multifiltrations satisfy the following crucial property.
\begin{theorem}
    \label{thm:intro-main}
    For any ordered set $Q$ and any $Q$-indexed filtration $X = (X^q)_q$ of simplicial complexes, there exists a $Q$-indexed filtration $T = (T^q)_q$ of spanning trees of $X$ satisfying $X^q \setminus T^q \subseteq X^{q'} \setminus T^{q'}$ for any pair $q \leq q'$ in $Q$.
\end{theorem}

This existence result for a cofiltration of spanning trees $T$ associated to $X$ allows us to associate to every edge $\sigma$ contained in $X$ but not in $T$ an upper set decomposable module $M_\sigma$.
We show that the direct sum $\bigoplus_\sigma M_\sigma$, where $\sigma$ ranges over any edge $\sigma$ in $X$ but not in $T$, admits a canonical epimorphism of persistence modules onto $Z_1(X)$, and subsequently onto $H_1(X)$.
\begin{theorem}
    Let $Q$ be any finite ordered set and $X = (X^q)_{q \in Q}$ a $Q$-indexed filtration of finite simplicial complexes.
    Then, given a choice of simplicial order on $X$ there exists an epimorphism $\varphi_1\colon M \to Z_1(X)$, where $M$ is an upper set decomposable module, constructed explicitly from the cofiltration of order-minimal spanning trees.
\end{theorem}

As a last step and in order to generalize that technique to higher homology groups $H_n(X)$, $n \in \mathbb N$, we introduce a direct generalization of spanning trees to higher dimensions.
We call this generalized structure an \emph{$n$-spanning complex} and show the existence by an application of Zorn's lemma.

\begin{theorem}
    For any simplicial complex~$X$ and any $n \in \mathbb N$ there exists a subcomplex $A$ such that $Z_n(A) = 0$ and $B_{n-1}(A) = B_{n-1}(X)$, i.e.\ $X$ admits an $n$-spanning complex.
\end{theorem}

\subsection{Organization}
In Section~\ref{sect:prelim}, we recall the basic definitions of partially ordered sets, simplicial complexes, simplicial homology, persistence modules, and persistent simplicial homology.
In Section~\ref{sect:span-trees}, we define our notion of spanning trees in terms of the boundaries and chains of subcomplexes.
Subsequently, we state an edge-exchange lemma (Lemma~\ref{lem:edge-exchange}) and the existence of a classical isomorphism associated to spanning trees (Proposition~\ref{prop:st-main}).

In Section~\ref{sect:cofil-st}, we investigate the previously introduced notion of spanning trees on their behavior under persistence.
In order to describe this behavior, we introduce the notion of a cofiltration of spanning trees (Definition~\ref{dfn:cofst}) and prove that they always exist for finite ordered sets and (globally) finite filtrations (Theorem~\ref{thm:main}).

The Section~\ref{sect:r-spans} is used to collect a variety of algebraic statements concerning so-called $k$-spans of set-valued diagrams over any small category.
These results are meant to be used in Section~\ref{sect:gen-h1}, where we construct from a cofiltration of spanning trees an upper set decomposable persistence module that enjoys an epimorphism onto the first persistent homology group (Theorem~\ref{thm:gen-h1}).

In the appendix, Section~\ref{sect:n-span}, we state the existence of $n$-spanning complexes, which are the natural generalization of the concept of spanning trees to higher dimensions.

\section{Preliminaries}%
\label{sect:prelim}

Throughout this paper, let $k$ be a principal ideal domain.

\begin{definition}
    Let $Q$ be a non-empty set and $R \subseteq Q \times Q$.
    \begin{enumerate}
        \item The relation~$R$ is \emph{reflexive} if $(q, q) \in R$ for all $q \in Q$.
        \item The relation~$R$ is \emph{transitive} if $(q, q'), (q', q'') \in R$ implies $(q, q'') \in R$.
        \item The relation~$R$ is \emph{antisymmetric} if $(q, q') \in R$ and $(q', q) \in R$ implies $q = q'$ for all $q, q' \in Q$.
        \item If $R$ is reflexive, transitive, and antisymmetric, then it is called a \emph{(partial) order} on $Q$. The pair $(Q, R)$ is called an \emph{ordered set}.
        \item If $R$ is an order and satisfies $(q, q') \in R$ or $(q', q) \in R$ for all $q, q' \in Q$, then $R$ is called a \emph{total order}. The pair $(Q, R)$ is called an \emph{totally ordered set}.
        \item If $R$ is an order, it is a \emph{well-ordering} if it is total and every subset of $Q$ has a minimal element w.r.t.\ $R$.
    \end{enumerate}
\end{definition}

\subsection{Simplicial complexes and simplicial homology}

\begin{definition}
    Let $V$ be a non-empty set.
    \begin{enumerate}
        \item An \emph{(abstract) simplicial complex} over $V$ is a set~$X$ of finite subsets of~$V$ such that
        \[ \sigma \in X \text{ and } \tau \subseteq \sigma \implies \tau \in X \]
        for all $\sigma, \tau \subset V$.
        \item The elements of a simplicial complex are called \emph{simplices}.
        A simplex of cardinality $n+1$ is called \emph{$n$-simplex}.
        \item The \emph{$n$-skeleton} of a simplicial complex~$X$ is defined by $\skel_n X = \bigcup_{j \leq n} X_j$.
        \item If $X$ is a simplicial complex over $V$, the set $V$ is denoted by $V(X)$ and is called the \emph{vertex set of $X$}.
        \item Let $X$ and $Y$ be simplicial complexes.
        A map $f\colon V(X) \to V(Y)$ is called \emph{simplicial map} if $f(\sigma) \in Y$ for all $\sigma \in X$.
    \end{enumerate}
\end{definition}

\begin{remark}
    Simplicial complexes with simplicial maps constitute the category $\Simp$ of simplicial complexes.
\end{remark}

For the purposes of (simplicial) homology, the vertex set $V(X)$ of a simplicial complex $X$ is endowed with a fixed total order $\preceq$.
We call such total order a \emph{vertex order}.
We assume that every simplicial complex has a vertex order.

Given an $n$-simplex $\sigma = \{ v_0, \dotsc, v_n \}$, we write $[v_0, \dotsc, v_n]$ for the simplex if $v_0 \preceq \dotsb \preceq v_n$.
Further, for $0 \leq i \leq n$, let $[v_0, \dotsc, \widehat{v_i}, \dotsc, v_n]$ denote the $(n-1)$-simplex that we obtain by removing the vertex $v_i$.

\begin{definition}%
    \label{defn:simp-hom}
    Let $X$ be a simplicial complex with a vertex order.
    \begin{enumerate}
        \item The \emph{chain module} is the free $k$-module~$C_n(X)$ generated by the $n$-simplices of $X$.
        \item The \emph{boundary homomorphism} of $X$ is defined by
        \begin{align*}
            \partial_n\colon C_n(X) &\longrightarrow C_{n-1}(X) \text, \\
            [v_0, \dotsc, v_n] &\longmapsto \sum_{i = 0}^n (-1)^n [v_0, \dotsc, \widehat{v_i}, \dotsc, v_n] \text.
        \end{align*}
        \item The \emph{$n$-cycle module} of $X$ is defined by $Z_n(X) = \Ker \partial_n$, and the \emph{$n$-boundary module} of $X$ is defined by $B_n(X) = \operatorname{Im} \partial_{n+1}$.
        \item The \emph{$n$-th homology module} of $X$ is the quotient $H_n(X) = Z_n(X)/B_n(X)$.
        \item If $f\colon X \to Y$ is a simplicial map, then the \emph{induced homomorphism} of $f$ is
        \[ H_n(f) = f_\ast\colon H_n(X) \to H_n(Y) \text, \quad \sigma + B_n(X) \mapsto f(\sigma) + B_n(Y) \text. \]
    \end{enumerate}
\end{definition}

\begin{remark}\strut
    \begin{enumerate}
        \item The $n$-th homology defines a functor from the category of simplicial complexes to the category of $k$-modules.
        \item Because $k$ is a principal ideal domain, the submodules $Z_n(X)$ and $B_n(X)$ are free.
    \end{enumerate}
\end{remark}

\subsection{Ordered simplicial complexes}

We extend the idea of endowing a vertex set with a total order to the whole simplicial complex.
Let $X$ be a simplicial complex and $\preceq$ be an order on the simplices of $X$.
The pair $(X, \mathord\preceq)$ is an \emph{ordered simplicial complex}.

Note that an ordered simplicial complex $(X, \mathord\preceq)$ comes with a natural vertex order by the correspondence
\[ V(X) \leftrightarrow \{ \{ v \} \mid v \in V \} \text. \]

Let $(X, \mathord\preceq_X)$ and $(Y, \mathord\preceq_Y)$ be ordered simplicial complexes.
A simplicial map $f\colon X \to Y$ is said to be \emph{ordered} or \emph{order-preserving} if
\[ \sigma \preceq_X \tau \implies f(\sigma) \preceq_Y f(\tau) \]
for all $\sigma, \tau \in X$.

\begin{remark}
    Ordered simplicial complexes with order-preserving simplicial maps constitute the category of ordered simplicial complexes.
    We denote this category by $\Simp_\ord$.
\end{remark}

\subsection{Persistence modules and persistent homology}

We give the definitions of persistence modules and of the persistent version of the simplicial homology theory of Definition~\ref{defn:simp-hom}.

\begin{definition}
    Let $(Q, \mathord\leq)$ be an ordered set.
    \begin{enumerate}
        \item A \emph{persistence module}~$M$ over $Q$ consists of a family $(M^q)_{q \in Q}$ of $k$-modules and a family $(\varphi^M_{q'q})_{q \leq q'}$ of $k$-homomorphisms $\varphi^M_{q'q}\colon M^q \to M^{q'}$ such that
        \[ \varphi^M_{q''q} = \varphi^M_{q''q'} \circ \varphi^M_{q'q} \]
        for all $q \leq q' \leq q''$ in $Q$.
        The homomorphisms $\varphi^M_{q'q}$ are called \emph{structure maps} of $M$.
        \item A \emph{persistent homomorphism} $f\colon M \to N$ is a family $(f^q)_{q \in Q}$ of $k$-ho\-mo\-mor\-phisms $f^q\colon M^q \to N^q$ such that
        \[ \varphi^N_{q'q} \circ f^q = f^{q'} \circ \varphi^M_{q'q} \]
        for all $q \leq q'$ in $Q$.
        \item A family $(X^q)_{q \in Q}$ of simplicial complexes over a common vertex set is a \emph{filtration} if $X^q \subseteq X^{q'}$ for all $q \leq q'$ in $Q$.
        \item Let $X = (X^q)_{q \in Q}$ be a filtration.
        For $n \in \mathbb N$ the \emph{$n$-dim.\ persistent homology module} of~$X$ is the persistence module~$M$ with
        \[ M^q = H_n(X^q) \]
        and structure maps defined by
        \[ \varphi^M_{q'q}\colon H_n(X^q) \longrightarrow H_n(X^{q'}) \text, \quad \sigma + B_n(X^q) \longmapsto \sigma + B_n(X^{q'}) \text, \]
        for $q \le q' \in Q$.
        This persistence module is denoted by $H_n(X)$.
    \end{enumerate}
\end{definition}

Furthermore, we call a set-valued functor on $Q$ a \emph{persistent set} and a corresponding natural transformation a \emph{persistent map}.
The category of persistent sets is denoted by $\Set^Q$.

\section{Spanning trees of simplicial complexes}%
\label{sect:span-trees}

Spanning trees of graphs are subgraphs that reflect the connectivity properties of the original graph while preserving acyclicity and means to identify cyclic paths in a graph \cite{BondyMurty2008,Kozlov2020,Lovász1977}.
As such, they are classically characterized by combinatiorial conditions.
We modify the standard definition of a spanning tree to facilitate the need of having an algebraic characterization in the following way:
\begin{itemize}
    \item We characterize spanning trees in terms of their cycle and boundary modules of simplicial homology, see Definition~\ref{def:st} below.
    \item The zeroth homology module of a spanning tree does not need to trivialize, i.e.\ a spanning tree is not necessarily connected.
\end{itemize}

Throughout this section, let $X$ be a simplicial complex with vertex set $V$.

\subsection{Homological spanning trees}

\begin{definition}[Homological spanning tree]%
    \label{def:st}
    A \emph{(homological) spanning tree} of~$X$ is a simplicial subcomplex $T$ of $X$ such that
    \begin{enumerate}[label=(\roman*)]
        \item\label{def:st:c1} $X_0 = T_0$,
        \item\label{def:st:c2} the boundary homomorphism $\partial_1\colon C_1(T) \to C_0(T)$ is injective, and
        \item\label{def:st:c3} $B_0(X) = B_0(T)$.
    \end{enumerate}
    The simplicial complex~$X$ is said to be \emph{spanned by $T$}.
\end{definition}

Condition~\ref{def:st:c1} states the natural assumption that any spanning tree of $X$ must have the same vertices (in the sense of 0-simplices) as $X$.
Further, condition~\ref{def:st:c2} requires that any spanning tree of $X$ is a graph-theoretic tree, i.e.\ the 1-dimensional cycle module of any spanning tree is trivial.

\begin{remark}
    Our main focus will be on condition~\ref{def:st:c3}, which states that a spanning tree spans the space $X$ in the following sense:
    Suppose that $T$ is a spanning tree of $X$ and let $u, v$ be vertices of $X$.
    Further, suppose that there exists a path in $X$ from $u$ to $v$, i.e.\ $u$ and $v$ lie in the same path component of $X$, or equivalent there is a 1-chain $c \in C_1(X)$ such that $\partial_1 c = u - v$.
    Then, the third condition implies that there exists a 1-chain $c' \in C_1(T)$ such that $\partial_1 c = \partial_1 c' \in B_0(T)$.
    Thus, both vertices $u$ and $v$ also lie in the same path component of $T$.

    More generally, if $c$ is a 1-chain in $X$, then there always exists a 1-chain~$c'$ in $T$ such that $c - c'$ is a cycle.
    It follows from Definition~\ref{def:st} \ref{def:st:c2} that the chain $c' \in C_1(T)$ is uniquely determined by $c \in C_1(X)$.
\end{remark}

\subsection{Operations on spanning trees}

Let $T$ be a spanning tree of $X$.
We pick some edge (1-simplex)~$e$ in $X$ that is not in $T$, and we ask whether there exists another spanning tree of $X$ that actually contains $e$.
It is clear that $A = T \cup \{ e \}$ still satisfies $B_0(X) = B_0(A)$ and $V(X) = V(A)$, and because of that the subcomplex $A$ still satisfies two out of three conditions for a spanning tree, see Definition~\ref{def:st} \ref{def:st:c1} and \ref{def:st:c3}.
However, it is also easy to see that $Z_1(A) \neq 0$.
Thus, $A$ does not satisfy the acyclicity condition of spanning trees.

Therefore, the actual question is whether we can find some edge $e'$ in $T$ such that $T \cup \{ e \} \setminus \{ e' \}$ is again a spanning tree.
We call the following lemma the \emph{edge-exchange lemma}.
A similar statement about spanning trees is known as the \emph{tree exchange property}, see \cite[p.\ 113, Exercise 4.3.2]{BondyMurty2008}.

\begin{lemma}%
    \label{lem:edge-exchange}
    Let $T$ be a spanning tree of the simplicial complex~$X$.
    For every edge~$\sigma \in X$ that is not contained in $T$, there exists an edge $\tau \in T$ such that $T \cup \{ \sigma \} \setminus \{ \tau \}$ is a spanning tree of~$X$.

    Moreover, if $c$ is a 1-chain in $T$ such that $\partial_1 c = \partial_1 \sigma$, then every simplex in $c$ is admissible as $\tau$.

    \begin{proof}
        Let $\sigma$ be an edge contained in $X$ but not in $T$.
        Since $T$ is a spanning tree, we find a 1-chain $c$ in $T$ such that $\partial_1 c = \partial_1 \sigma$.
        Moreover, we can choose $\lambda_0, \dotsc, \lambda_m \in \{ 1, -1 \}$ and distinct edges $\tau_0, \dotsc, \tau_m \in T$ such that $c = \sum_{i = 0}^m \lambda_i \tau_i$.
        We set $\tau = \tau_0$, i.e.\ $c = \lambda_0 \tau + \sum_{i = 1}^m \lambda_i \tau_i$.

        Let $b \in B_0(X)$.
        Then, we can choosse a 1-chain $c' = \mu_0 \tau + \sum_{i = 1}^l \mu_i \sigma_i$ in $T$ such that $\partial_1 c' = b$, where the simplices $\tau, \sigma_1, \dotsc, \sigma_m$ are distinct.
        By replacing $\tau$ in $c'$ with the 1-chain $\sum_{i = 1}^m \lambda_i \tau_i - \sigma$, which is contained in $X \setminus T$, we obtain the chain
        \[ c'' = \sum_{i = 1}^m \mu_0 \lambda_i \tau_i - \mu_0 \sigma + \sum_{i = 1}^m \mu_i \sigma_i \text. \]
        Since the chains $\sum_{i = 1}^m \lambda_i \tau_i$ and $\sigma$ have the same boundary, the boundary of $c''$ coincides with the boundary of $c'$.
        Further, the chain $c''$ is, by construction, contained in $T'$.
        Hence, $\partial_1 c'' = \partial_1 c' = b \in B_0(T')$.
        Since $T' \subseteq X$, this proves in particular $B_0(X) = B_0(T')$ and $B_0(T) = B_0(T')$.

        Let $z \in Z_1(T + \sigma - \tau)$ and choose $\mu_0, \mu_1, \dotsc, \mu_l \in k$ and edges $\tau_1, \dotsc, \tau_l$ distinct from $\tau$ in $T$ such that $z = \mu_0 \sigma + \sum_{j = 1}^l \mu_j \tau_j$.
        Then, we have $\mu_0 \neq 0$ unless $z = 0$.
        Therefore, we assume that $z \neq 0$.

        Then, we have
        \[ \mu_0 (u - v) = \mu_0 \partial_1 \sigma = -\partial_1 \sum_{j = 1}^l \mu_j \tau_j  \]
        for vertices distinct $u$ and $v$ in $X$.
        Furthermore, $\mu_0 \neq 0$ is invertible since there exists some chain $c \in C_1(T)$ satisfying $\partial_1 c = u - v$ and $\sum_{j = 1}^l \mu_j \tau_j = \mu_0 c$ by the injectivity of $\partial_1$ on $C_1(T)$.
        Thus, $u - v = -\partial_1 \sum_{j = 1}^l \mu_j \mu_0^{-1} \tau_j$.

        However, this means that some of the $\tau_1, \dotsc, \tau_j$ must coincide with $\tau$ by construction.
        This is a contradiction to the assumption that $\tau_1, \dotsc, \tau_j$ are distinct from $\tau$.
        Hence, $Z_1(T + \sigma - \tau) = 0$.
    \end{proof}
\end{lemma}

\begin{definition}
    Let $T$ be a spanning tree of $X$.
    A pair $(\sigma, \tau)$ of edges $\sigma \in X \setminus T$ and $\tau \in T$ is called an \emph{edge exchange pair} of $T$ in $X$ if $T \cup \{ \sigma \} \setminus \{ \tau \}$ is a spanning tree.

    If $(\sigma, \tau)$ is an edge-exchange pair, then $T \cup \{ \sigma \} \setminus \{ \tau \}$ is denoted by $T + \sigma - \tau$.
\end{definition}

\subsection{Homology of spanning trees}

Let $X$ be a simplicial complex and $T$ a spanning tree of $X$.
The relative homology of the pair $(X, T)$ is intimately connected to the homology of~$X$.
This is a standard result in combinatorial algebraic topology, and proofs of such results can be found in  \cite[Theorem~1.9.4]{Diestel2017}, \cite[Theorem~2.33]{Kozlov2020}, and \cite[Theorem~3, p.\ 140]{Spanier1982}.

\begin{proposition}%
    \label{prop:st-main}
    Let $T$ be a spanning tree of $X$.
    The homomorphism
    \[ j\colon Z_1(X) \longrightarrow C_1(X, T) \text, \quad z \longmapsto z + C_1(T) \text, \]
    is an isomorphism of $k$-modules.

    \begin{proof}
        Let $z \in Z_1(X)$ such that $z \in C_1(T)$.
        Then, $z = 0$ by $Z_1(T) = 0$.
        Hence, $j$ is a monomorphism.

        In order to prove the surjectivity of $j$, let $c \in C_1(T)$.
        By $B_0(T) = B_0(X)$, we find $c' \in C_1(T)$ such that $\partial_1 c = \partial_1 c'$, that is~$c - c' \in Z_1(X)$.
        Then, we have
        \[ j(c - c') = c - c' + C_1(T) = c + C_1(T) \]
        since $c' \in C_1(T)$.
        Thus, the homomorphism~$j$ is surjective, and it follows $Z_1(X) \cong C_1(X, T)$.
    \end{proof}
\end{proposition}

\subsection{Minimality of spanning trees}

We assume that every ordered simplicial complex $(X, \mathord\preceq)$ satisfies the following properties:
\begin{enumerate}[label=(\roman*)]
    \item The total order~$\mathord\preceq$ is a well-order, i.e.\ every subset $A$ of $X$ admits a minimal element.
    \item If $\sigma$ is a simplex of $X$ and $\tau$ a face of $\sigma$, then $\tau \preceq \sigma$.
\end{enumerate}

It turns out that every finite, ordered simplicial complex admits a special, distinguished spanning tree that is minimal in a certain sense w.r.t.\ the well-ordering~$\preceq$.
We characterize these spanning trees using lexicographic orders.

\begin{definition}[Lexicographic orders]
    Let $(X, \mathord\preceq)$ be a well-ordered set, e.g.\ $(X, \mathord\preceq)$ is an ordered simplicial complex.
    The \emph{lexicographic order}~$\mathord\preceq_\lex$ on $2^X$ is given by
    \[ A \preceq_\lex B \iff A = B \text{ or } \min((A \setminus B) \cup (B \setminus A)) \in A \]
    for subsets $A, B \subseteq X$.
\end{definition}

\begin{remark}
    Let $A, B \subseteq X$, where $X$ is an arbitrary well-ordered set~$X$.
    Then, $A \preceq_\lex B$ if and only if there exists $x \in A \setminus B$ such that for every $y \in (A \setminus B) \cup (B \setminus A)$ it holds $x \preceq y$.
\end{remark}

\begin{lemma}%
    \label{lem:lex-is-total}
    The lexicographic order~$\preceq_\lex$ is a total order.
\end{lemma}

\begin{remark}
    If $X$ is a finite set, its power set~$2^X$ is also finite, and then the lexicographic order~$\preceq_\lex$ is a well-ordering of $2^X$.
\end{remark}

\begin{example}
    Consider the totally ordered set $X = \{ 1 < 2 < 3 \}$ and endow its power set
    \[ 2^X = \big\{ \emptyset, \{1\}, \{2\}, \{3\}, \{1, 2\}, \{1,3\}, \{2, 3\}, \{1,2,3\} \big\} \]
    with the lexicographic order $<_\lex$.
    Then it holds $\{ 1, 2 \} <_\lex \{ 1, 3 \} <_\lex \{ 2, 3 \}$.
    This motivates the term \enquote{lexicographic order} for the order defined above.

    See Figure~\ref{fig:123-lex-order} for a table of all relations in $2^X$.
    We note that $A \subseteq B$ implies $B <_\lex A$ and therefore the element $\emptyset$ is maximal. Analogously, $X \in 2^X$ is minimal w.r.t.\ the lexicographic order.
\end{example}

\begin{figure}
    \centering
    \begin{tabular}{r|cccccccc}
                    & $\emptyset$ & $1$    & $2$    & $3$    & $12$   & $13$   & $23$   & $123$  \\ \hline
        $\emptyset$ & $=$         & $\geq$ & $\geq$ & $\geq$ & $\geq$ & $\geq$ & $\geq$ & $\geq$ \\
        $1$         &             & $=$    & $\leq$ & $\leq$ & $\geq$ & $\geq$ & $\leq$ & $\geq$ \\
        $2$         &             &        & $=$    & $\leq$ & $\geq$ & $\geq$ & $\geq$ & $\geq$ \\
        $3$         &             &        &        & $=$    & $\geq$ & $\geq$ & $\geq$ & $\geq$ \\
        $12$        &             &        &        &        & $=$    & $\leq$ & $\leq$ & $\geq$ \\
        $13$        &             &        &        &        &        & $=$    & $\leq$ & $\geq$ \\
        $23$        &             &        &        &        &        &        & $=$    & $\geq$ \\
        $123$       &             &        &        &        &        &        &        & $=$
    \end{tabular}

    \caption{The relations of the lexicographic order on $2^X$ for $X = \{ 1 < 2 < 3 \}$, where the entry \enquote{$\geq$} denotes \enquote{$\text{row label} \geq \text{column label}$} and the entry \enquote{$\leq$} denotes \enquote{$\text{row label} \leq \text{column label}$}.}
    \label{fig:123-lex-order}
\end{figure}

Let $(X, \mathord\preceq)$ be an finite, ordered simplicial complex.
Define the set
\[ \mathcal T(X) = \{ T \mid \text{$T$ is a spanning tree of $X$} \} \text. \]
As a subset of the power set of $2^X$, we endow it with the lexicographic ordering~$\preceq_\lex$.
Since $X$ is a finite simplicial complex, the set $\mathcal T(X)$ is also a finite, total ordered set.
Therefore, $\mathcal T(X)$ admits a unique, minimal element.

\begin{definition}[Order-minimal spanning tree]
    The minimal element of $\mathcal T(X)$ is the \emph{order-minimal spanning tree} of~$(X, \mathord\preceq)$.
\end{definition}

\begin{example}
    Consider the triangle complex $X = \{ \emptyset, 1, 2, 3, 12, 13, 23 \}$.
    This complex admits three spanning trees $T_1$, $T_2$, and $T_3$, which are depicted in Figure~\ref{fig:triangle-st}.
    When we endow the triangle complex with the lexicographic order induced by $1 \preceq 2 \preceq 3$, we see that $T_1$ is the order-minimal spanning tree of $X$.
\end{example}

\begin{figure}
    \centering
    \begin{tikzpicture}[
        dot/.style={radius=2.5pt, draw=white, fill=black}
    ]
        \draw[dashed]
            (0, 0) -- (1, 0)
            (4, 0) -- (3.5, 1)
            (6, 0) -- (6.5, 1);
        \draw
            (0, 0) -- (0.5, 1) -- (1, 0)
            (3.5, 1) -- (3, 0) -- (4, 0)
            (6.5, 1) -- (7, 0) -- (6, 0)
            (9, 0) -- (9.5, 1) -- (10, 0) -- cycle;
        \draw[dot]
            (0, 0) circle {} node[anchor=south east] {$v_1$}
            (1, 0) circle {} node[anchor=south west] {$v_2$}
            (3, 0) circle {} node[anchor=south east] {$v_1$}
            (4, 0) circle {} node[anchor=south west] {$v_2$}
            (6, 0) circle {} node[anchor=south east] {$v_1$}
            (7, 0) circle {} node[anchor=south west] {$v_2$}
            (9, 0) circle {} node[anchor=south east] {$v_1$}
            (10, 0) circle {} node[anchor=south west] {$v_2$}
            (0.5, 1) circle {} node[anchor=south west] {$v_3$}
            (3.5, 1) circle {} node[anchor=south west] {$v_3$}
            (6.5, 1) circle {} node[anchor=south west] {$v_3$}
            (9.5, 1) circle {} node[anchor=south west] {$v_3$};
        \node at (0.5, -0.5) {$T_1$};
        \node at (3.5, -0.5) {$T_2$};
        \node at (6.5, -0.5) {$T_3$};
        \node at (9.5, -0.5) {$F_3$};
    \end{tikzpicture}

    \caption{The spanning trees $T_1$, $T_2$, and $T_3$ of the triangle complex $F_3$. The edges not contained in the spanning trees are indicated as dashed lines.}
    \label{fig:triangle-st}
\end{figure}

\begin{lemma}%
    \label{lem:min-edge-exchange}
    Let $T$ be the order-minimal spanning tree of~$(X, \mathord\preceq)$.
    For every edge exchange pair $(\sigma, \tau)$ of $T$ it holds $\tau \preceq \sigma$.

    \begin{proof}
        Suppose there exists an edge exchange pair~$(\sigma, \tau)$, $\sigma \in X \setminus T$ and $\tau \in T$, such that $\sigma \prec \tau$.
        Then, for the spanning tree~$T' = T + \sigma - \tau$ we have
        \[ A = (T' \setminus T) \cup (T \setminus T') = \{ \sigma, \tau \} \text. \]
        Thus, by assumption, the minimal element of~$A$ is~$\sigma$, and $\sigma$ is contained in $T'$.
        This implies $T' \prec_\lex T$.
        Hence, $T$ is not order-minimal.
    \end{proof}
\end{lemma}

\begin{proposition}%
    \label{prop:minimal-charact}
    A spanning tree~$T$ of~$X$ is order-minimal w.r.t.\ $\preceq$ if and only if there is no edge-exchange pair $(\sigma, \tau)$ such that $\sigma \preceq \tau$.

    \begin{proof}
        The assertion follows from Lemma~\ref{lem:min-edge-exchange} together with Lemma~\ref{lem:edge-exchange}:
        If $T$ is order-minimal, every edge exchange pair $(\sigma, \tau)$ must satisfy $\sigma \preceq \tau$.
        Vice-versa, if $T$ is not order-minimal, let $T'$ be an order-minimal spanning tree of $X$ and pick an edge-exchange pair $(\sigma, \tau)$ where $\sigma$ is an edge in $T'$ but not in $T$.
    \end{proof}
\end{proposition}

\begin{remark}
    The characterization of order-minimal spanning trees in terms of their edge-exchange pairs in conjunction with the edge-exchange lemma (Lemma~\ref{lem:edge-exchange}) provides a tactic to reduce a given spanning tree to the order-minimal one.

    Let $T$ be any spanning tree of~$X$.
    Then, if $T$ is not order-minimal, we can iteratively pick an edge exchange pair $(\sigma, \tau)$ such that $\sigma \preceq \tau$, and replace $T$ by $T + \sigma - \tau$.

    This method iteratively reduces the spanning tree~$T$ until $T$ is order-minimal by the characterization above.
\end{remark}

\begin{definition}[Leading simplex]
    Let $n \in \mathbb N$ and let $c = \lambda_1 \sigma_1 + \dotsb + \lambda_m \sigma_m$ be a non-zero $n$-chain of $X$.
    Suppose that $\sigma_1, \dotsc, \sigma_m$ are distinct, $\lambda_1, \dotsc, \lambda_m$ are non-zero, and $\sigma_1 \preceq \dotsb \preceq \sigma_m$.
    The \emph{leading simplex}~$\LS(c) = \LS_{\mathord\preceq}(c)$ of $c$ is $\sigma_m$.
\end{definition}

\begin{proposition}
    \label{prop:ls-complement}
    Let $T$ be an order-minimal spanning tree of $(X, \mathord\preceq)$.
    For every non-zero $z \in Z_1(X)$, it holds $\LS(z) \in X \setminus T$.

    \begin{proof}
        Let $z \in Z_1(X)$.
        Then, there exist simplices $\sigma_1, \dotsc, \sigma_m$ in $X$ but not in $T$ and non-zero $\lambda_1, \dotsc, \lambda_m \in k$ such that
        \[ z = \sigma_1 \lambda_1 + \dotsc + \lambda_m \sigma_m + c \]
        for some 1-chain $c$ in $T$.
        We choose for every $i \in \{ 1, \dotsc, m \}$ a 1-chain $c_i$ in $T$ such that $\partial_1 c_i = \partial_1 \sigma_i$.
        Then, we have
        \[ z = \lambda_1 (\sigma_1 - c_1) + \dotsc + \lambda_m (\sigma_m - c_m) \in Z_1(X) \text. \]

        Since we can construct an edge-exchange pair $(\sigma_i, \tau_i)$ of $T$ for any $i \in \{ 1, \dotsc, m \}$, it follows $\sigma_i \preceq \tau_i$.
        In particular, every $1$-simplex in $c_i$ is admissible as $\tau_i$, cf.\ Lemma~\ref{lem:edge-exchange}.
        Therefore, $\LS(\sigma_i - c_i) = \sigma_i$.
        Hence, we have
        \begin{align*}
            \LS(z) = \LS(\lambda_1 (\sigma_1 - c_1) + \dotsc + \lambda_m (\sigma_m - c_m)) &\in \LS \{ \sigma_i - c_i \mid i \in \{ 1, \dotsc, m \} \} \\
            &= \{ \sigma_1, \dotsc, \sigma_m \} \subseteq X \setminus T \text.
            \qedhere
        \end{align*}
    \end{proof}
\end{proposition}

\section{Filtrations and cofiltrations of spanning trees}
\label{sect:cofil-st}

Let $Q = (Q, \mathord\leq)$ be a finite ordered set and $X = (X^q)_{q \in Q}$ a $Q$-indexed filtration of finite simplicial complexes.
In this section, we answer the question of whether there exists a persistent analogue of spanning trees for filtrations of simplicial complexes.

\subsection{Filtrations}
First, we assume that there exists a filtration $T = (T^q)_{q \in Q}$ such that $T^q$ is a spanning tree of $X^q$ for every $q \in Q$.
We say that $T$ is a \emph{filtration of spanning trees} of $X$.

We consider the persistence module $C_1(X, T)\colon Q \to \Mod{k}$ with $C_1(X, T)(q) = C_1(X^q, T^q)$ for $q \in Q$ and structure maps
\[ C_1(X^q, T^q) \to C_1(X^{q'}, T^{q'}) \text, \quad c + C_1(T^q) \mapsto c + C_1(T^{q'}) \text. \]

\begin{lemma}%
    \label{lem:rel-hom-wf}
    The structure maps of $C_1(X, T)$ are well-defined module homomorphisms.

    \begin{proof}
        Let $q \leq q' \in Q$ be a pair.
        The well-definedness of the structure map from $C_1(X^q, T^q)$ to $C_1(X^{q'}, T^{q'})$ follows immediately from $T^q \subseteq T^{q'}$.

        Let $c \in C_1(X^q)$ such that $c$ is also a chain in $T^{q'}$.
        Then, we choose non-zero $\lambda_1, \dotsc, \lambda_m \in m$ and distinct simplices $\sigma_1, \dotsc, \sigma_m$ in $T^{q'} \cap X^{q}$ such that $c = \sum_{i = 1}^m \lambda_i \sigma_i$.
        We suppose without loss of generality and to the contrary that $\sigma_1 \notin T^q$.
        Then, $T^q \cup \{ \sigma_1 \}$ is a subcomplex of $X^q$ with non-trivial cycle module.
        However, $T^q \cup \{ \sigma_1 \}$ is also a subcomplex of $T^{q'}$ since $\sigma_1 \in T^{q'}$.
        Thus, we have a contradiction since $T^{q'}$ has trivial cycle module.
        It follows that $\sigma_1, \dotsc, \sigma_m$ are contained in $T^q$.
        Therefore, $c$ is a 1-chain in $T^q$, and the kernel of the structure map from $C_1(X^q, T^q)$ to $C_1(X^{q'}, T^{q'})$ is trivial.
    \end{proof}
\end{lemma}

Now, we have for any $q \in Q$ an $k$-module isomorphism $j^q\colon Z_1(X^q) \to C_1(X^q, T^q)$.
These degreewise isomorphisms extend to a persistent isomorphism.

\begin{proposition}%
    \label{prop:pers-filt}
    Suppose there exists a filtration $T = (T^q)_{q \in Q}$ of simplicial complexes such that $T^q$ is a spanning tree of $X^q$ for each $q \in Q$.
    The persistent homomorphism $j\colon Z_1(X) \to C_1(X, T)$ that is defined at $q \in Q$ by
    \[ j^q\colon Z_1(X^q) \to C_1(X^q, T^q) \text, \qquad z \mapsto z + C_1(T^q) \text, \]
    is an isomorphism of persistence modules.

    \begin{proof}
        The assertion is equivalent to the commutativity of the diagram
        \[ \begin{tikzcd}[ampersand replacement=\&]
            Z_1(X^q) \& C_1(X^q, T^q) \\
            Z_1(X^{q'}) \& C_1(X^{q'}, T^{q'})
            \arrow[from=1-1, to=1-2, "j^q"]
            \arrow[from=1-1, to=2-1]
            \arrow[from=1-2, to=2-2]
            \arrow[from=2-1, to=2-2, "j^{q'}"']
        \end{tikzcd} \]
        where the structure map of $Z_1(X)$ from $Z_1(X^q)$ to $Z_1(X^{q'})$ is the canonical inclusion map.
        However, the commutativity of the diagram follows immediately.
    \end{proof}
\end{proposition}

\begin{remark}
    If $Q$ is totally ordered, we can always choose a subfiltration of spanning trees.
\end{remark}

\begin{figure}
    \begin{tikzpicture}[
        commutative diagrams/every diagram,
        anchor box/.style={minimum width=2.0cm, minimum height=1.8cm, anchor=south west, at={(-5mm, -4mm)}},
        bounding box/.style={},
        dot/.style={radius=2.5pt, draw=white, fill=black},
        hidden edge/.style={color=gray, dashed},
        edge/.style={draw=black},
        emph path/.style={
            start angle=90,
            end angle=400,
            radius=1.5mm,
            yshift=1.5mm,
            thick,
            arrows={->[length=0.7mm]},
            draw=orange
        },
        emph path label/.style={},
        baseline=(m-1-1.north)]
        \matrix[matrix anchor=west, column sep=5mm, row sep=5mm] at (0, 0) {
            \draw[hidden edge]
                (0, 0) -- (1, 0) -- (1, 1);
            \draw[edge] (0, 0) -- (0, 1) -- (1, 1) -- cycle;
            \path[dot, fill=blue] (0, 0) circle[] node[anchor=south east] {$v_1$};
            \path[dot, fill=red] (1, 0) circle[] node[anchor=south east] {$v_2$};
            \path[dot, fill=teal] (0, 1) circle[] node[anchor=south east] {$v_3$};
            \path[dot, fill=violet] (1, 1) circle[] node[anchor=south east] {$v_4$};
            \node[anchor box] (m-1-1) {};
            \&
            \draw[edge] (0, 0) -- (0, 1) -- (1, 1) -- (1, 0) -- cycle
                (0, 0) -- (1, 1);
            \path[dot, fill=blue] (0, 0) circle[] node[anchor=south east] {$v_1$};
            \path[dot, fill=red] (1, 0) circle[] node[anchor=south east] {$v_2$};
            \path[dot, fill=teal] (0, 1) circle[] node[anchor=south east] {$v_3$};
            \path[dot, fill=violet] (1, 1) circle[] node[anchor=south east] {$v_4$};
            \node[anchor box] (m-1-2) {};
            \&
            \draw[edge] (0, 0) -- (0, 1) -- (1, 1) -- (1, 0) -- cycle
                (0, 0) -- (1, 1);
            \path[dot, fill=blue] (0, 0) circle[] node[anchor=south east] {$v_1$};
            \path[dot, fill=red] (1, 0) circle[] node[anchor=south east] {$v_2$};
            \path[dot, fill=teal] (0, 1) circle[] node[anchor=south east] {$v_3$};
            \path[dot, fill=violet] (1, 1) circle[] node[anchor=south east] {$v_4$};
            \node[anchor box] (m-1-3) {};
            \\
            \draw[hidden edge]
                (0, 0) -- (1, 0) -- (1, 1) -- cycle;
            \draw[edge]
                (0, 0) -- (0, 1) -- (1, 1);
            \path[dot, fill=blue] (0, 0) circle[] node[anchor=south east] {$v_1$};
            \path[dot, fill=red] (1, 0) circle[] node[anchor=south east] {$v_2$};
            \path[dot, fill=teal] (0, 1) circle[] node[anchor=south east] {$v_3$};
            \path[dot, fill=violet] (1, 1) circle[] node[anchor=south east] {$v_4$};
            \node[anchor box] (m-2-1) {};
            \&
            \draw[hidden edge]
                (0, 0) -- (1, 1);
            \draw[edge]
                (0, 0) -- (0, 1) -- (1, 1) -- (1, 0) -- cycle;
            \path[dot, fill=blue] (0, 0) circle[] node[anchor=south east] {$v_1$};
            \path[dot, fill=red] (1, 0) circle[] node[anchor=south east] {$v_2$};
            \path[dot, fill=teal] (0, 1) circle[] node[anchor=south east] {$v_3$};
            \path[dot, fill=violet] (1, 1) circle[] node[anchor=south east] {$v_4$};
            \node[anchor box] (m-2-2) {};
            \&
            \draw[edge] (0, 0) -- (0, 1) -- (1, 1) -- (1, 0) -- cycle
                (0, 0) -- (1, 1);
            \path[dot, fill=blue] (0, 0) circle[] node[anchor=south east] {$v_1$};
            \path[dot, fill=red] (1, 0) circle[] node[anchor=south east] {$v_2$};
            \path[dot, fill=teal] (0, 1) circle[] node[anchor=south east] {$v_3$};
            \path[dot, fill=violet] (1, 1) circle[] node[anchor=south east] {$v_4$};
            \node[anchor box] (m-2-3) {};
            \\
            \draw[hidden edge]
                (0, 0) -- (0, 1) -- (1, 1) -- (1, 0) -- cycle
                (0, 0) -- (1, 1);
            \path[dot, fill=blue] (0, 0) circle[] node[anchor=south east] {$v_1$};
            \path[dot, fill=red] (1, 0) circle[] node[anchor=south east] {$v_2$};
            \path[dot, fill=teal] (0, 1) circle[] node[anchor=south east] {$v_3$};
            \path[dot, fill=violet] (1, 1) circle[] node[anchor=south east] {$v_4$};
            \node[anchor box] (m-3-1) {};
            \&
            \draw[hidden edge]
                (0, 0) -- (1, 1) -- (0, 1) -- cycle;
            \draw[edge] (0, 0) -- (1, 0) -- (1, 1);
            \path[dot, fill=blue] (0, 0) circle[] node[anchor=south east] {$v_1$};
            \path[dot, fill=red] (1, 0) circle[] node[anchor=south east] {$v_2$};
            \path[dot, fill=teal] (0, 1) circle[] node[anchor=south east] {$v_3$};
            \path[dot, fill=violet] (1, 1) circle[] node[anchor=south east] {$v_4$};
            \node[anchor box] (m-3-2) {};
            \&
            \draw[hidden edge]
                (0, 0) -- (0, 1) -- (1, 1);
            \draw[edge] (0, 0) -- (1, 0) -- (1, 1) -- cycle;
            \path[dot, fill=blue] (0, 0) circle[] node[anchor=south east] {$v_1$};
            \path[dot, fill=red] (1, 0) circle[] node[anchor=south east] {$v_2$};
            \path[dot, fill=teal] (0, 1) circle[] node[anchor=south east] {$v_3$};
            \path[dot, fill=violet] (1, 1) circle[] node[anchor=south east] {$v_4$};
            \node[anchor box] (m-3-3) {};
            \\
        };
        \path[commutative diagrams/.cd, every arrow, every label]
            (m-3-1) edge (m-3-2)
            (m-3-1) edge (m-2-1)
            (m-3-2) edge (m-3-3)
            (m-3-2) edge (m-2-2)
            (m-3-3) edge (m-2-3)
            (m-2-1) edge (m-2-2)
            (m-2-1) edge (m-1-1)
            (m-2-2) edge (m-2-3)
            (m-2-2) edge (m-1-2)
            (m-2-3) edge (m-1-3)
            (m-1-1) edge (m-1-2)
            (m-1-2) edge (m-1-3);
    \end{tikzpicture}

    \caption{A $3\times 3$-filtration of simplicial complexes, which does not admit a subfiltration of spanning trees.}
    \label{fig:no-subfil}
\end{figure}

\begin{example}
    \label{ex:no-subfil}
    If the ordered set $Q$ is not total, the existence of a subfiltration of spanning trees for any filtration over $Q$ is not guaranteed.
    A counterexample of such filtration over $\mathbb N^2$ is depicted in Figure~\ref{fig:no-subfil}.
    This failure of existence is also reflected in the fact that the decomposition structure of persistence modules over non-totally ordered sets is not trivial \cite{CarlssonZomorodian2009}.
\end{example}

\subsection{Cofiltrations}

\begin{definition}[Cofiltrations of spanning trees]
    \label{dfn:cofst}
    Let $(T^q)_{q \in Q}$ be a family of simplicial complexes such that $T^q$ is a spanning tree of $X^q$ for any $q \in Q$.
    The family $(T^q)_{q \in Q}$ is a \emph{cofiltration of spanning trees} of $X$ if $X^q \setminus T^q \subseteq X^{q'} \setminus T^{q'}$ set-theoretically for any pair $q \leq q'$ in $Q$.
\end{definition}

\begin{remark}
    A similar but more general notion of cofiltrations was introduced in \cite{PatelRask2022}:
    Given a finite simplicial complex $X$, a \emph{supcomplex}~$A$ of $K$ is a subset of $K$ closed under taking cofaces.
    Then, a \emph{cofiltration} over $K$ is a functor from $(Q, \mathord\leq)$ into the set of supcomplexes of $K$, ordered by inclusion.
    In this sense, every cofiltration of spanning trees is a cofiltration in the sense of Patel and Rask, cf.\ \cite{PatelRask2022}.
\end{remark}

\begin{lemma}%
    \label{lem:classic-comp-incl}
    Let $(X, \mathord\preceq)$ be an ordered simplicial complex and $A$ a subcomplex of $X$.
    Interpret $A$ as an ordered simplicial complex, with ordering given by the restriction of $\mathord\preceq$.
    For order-minimal spanning trees $T_A$ and $T_X$ of $A$ and $X$, respectively, it holds $A \setminus T_A \subseteq X \setminus T_X$ set-theoretically.

    \begin{proof}
        Let $e \in A \setminus T_A$.
        Then, we find $c = \sum_{i = 1}^m \lambda_i e_i \in C_1(T_A)$ such that the chain $z = c - e$ is a cycle in $A$.
        By Lemma~\ref{prop:ls-complement}, the leading simplex of $z$ must be contained in $A \setminus T_A$, and the only possibility is $\LS(z) = e$.
        Since $T_X$ is an order-minimal spanning tree of $X$, it further follows that $\LS(z) \in X \setminus T_X$.
        Hence, we have $e \in X \setminus T_X$.
    \end{proof}
\end{lemma}

Our main result in this section is that we can use the above lemma to obtain, for every filtration of finite simplicial complexes over a finite ordered set, a corresponding cofiltration of spanning trees, cf.\ Theorem~\ref{thm:intro-main}.

\begin{theorem}
    \label{thm:main}
    Let $X = (X^q)_{q \in Q}$ be a filtration of finite simplicial complexes.
    There exists a cofiltration of spanning trees $T$ of $X$, i.e.\ a family $T = (T^q)_{q \in Q}$ of simplicial complexes with $X^q \setminus T^q \subseteq X^{q'} \setminus T^{q'}$ for any pair $q \leq q' \in Q$.

    \begin{proof}
        We choose for any $q \in Q$ an order-minimal spanning tree $T^q$ of $X^q$, which exists since $X^q$ is finite.
        Let $q \leq q' \in Q$.
        We know by assumption that $X^q$ is a subcomplex of $X^{q'}$.
        Thus, by Lemma~\ref{lem:classic-comp-incl}, we have $X^q \setminus T^q \subseteq X^{q'} \setminus T^{q'}$.
        Hence, $(T^q)_{q \in Q}$ is a cofiltration of spanning trees of $(X^q)_{q \in Q}$.
    \end{proof}
\end{theorem}

\subsection{Functoriality of order-minimal spanning trees}

Lastly, we show that we can extend the above result in a functorial fashion for order- and dimension-preserving simplicial maps.
A morphism of simplicial complexes~$f \colon X \to Y$ is called \emph{dimension-preserving} if $\dim \sigma = \dim f(\sigma)$.
Let $\Simp_{\ord, \dim}$ denote the subcategory of ordered simplicial complexes with dimension- and order-preserving simplicial maps as morphisms.

\begin{definition}
    Let $A \subseteq X$ be a pair of simplicial complexes and $n \in \mathbb N$.
    The \emph{$n$-difference of $X$ and $A$} is the set
    \[ X \boxminus_n A = X_0 \cup \dotsc \cup X_{n-1} \cup (X_n \setminus A_n) \text. \]
\end{definition}

It is immediately clear that $X \boxminus_n A$ is again a simplicial complex contained in $X$.
Further, we note that the $(n-1)$-skeleton of $X \boxminus_n A$ coincides with the $(n-1)$-skeleton of $X$.

\begin{theorem}%
    \label{thm:comp-inc-fct}
    Let $f\colon (X, \mathord\leq_X) \to (Y, \mathord\leq_Y)$ be an order- and dimension-preserving simplicial map.
    Suppose that $T_X$ is an order-minimal spanning tree of $X$ and $T_Y$ is an order-minimal spanning tree of $Y$.
    Then, $f(X \boxminus_1 T_X) \subseteq Y \boxminus_1 T_Y$.

    Hence, $f$ restricts to a simplicial map $X \boxminus_1 T_X \to Y \boxminus_1 T_Y$.

    \begin{proof}
        Let $\lambda_1, \dotsc, \lambda_m \neq 0$ and $\sigma_1, \dotsc, \sigma_m$ be distinct $1$-simplices in $X$ satisfying $\sigma_1 \preceq_X \dotsb \preceq_X \sigma_m$.
        That is, for the chain $c = \lambda_1 \sigma_1 + \dotsb + \lambda_m \sigma_m$ we have $\LS(c) = \sigma_m$.

        Since $f$ is order-preserving, it follows $f(\sigma_1) \preceq_Y \dotsb \preceq_Y f(\sigma_m)$.
        Moreover, since $f$ is dimension-preserving, every simplex $f(\sigma_i)$ is an $1$-simplex for $i \in \{ 1, \dotsc, m \}$.

        Let $f_\sharp\colon C_1(X) \to C_1(Y)$ be the associated homomorphism of chain modules of~$f$.
        Therefore, we have $f_\sharp(c) = \lambda_1 f(\sigma_1) + \dotsb + \lambda_m f(\sigma_m)$.
        Hence, $\LS(f_\sharp(c)) = f(\sigma_m) = f(\LS(c))$.

        Now, let $\sigma$ be an $1$-simplex in $X$ but not in $T_X$.
        Then, we choose a chain $c$ in $T_X$ such that $c - \sigma \in Z_1(X)$.
        Since the spanning tree $T_X$ is order-minimal, we have $\LS(c - \sigma) \in X \setminus T_X$.
        Thus, it holds $\LS(c - \sigma) = \sigma$.

        By the considerations above, we have $\LS(f_\sharp(c - \sigma)) = f(\LS(c - \sigma)) = f(\sigma)$.
        Since $T_Y$ is also order-minimal and $f_\sharp(c - \sigma)$ is by definition a non-zero $1$-cycle in $Y$, it follows $\LS(f_\sharp(c - \sigma)) \in Y \setminus T_Y$.
        Hence, $f(\sigma) \in Y \boxminus_1 T_Y$, which proves the assertion.
    \end{proof}
\end{theorem}

\begin{remark}%
    \label{rem:mono-is-dim-preserv}
    Every monomorphism of simplicial complexes is necessarily dimension-preserving.
    Hence, if $f\colon (X, \mathord\leq_X) \to (Y, \mathord\leq_Y)$ is a monomorphism, it follows $f(X \setminus T_X) \subseteq Y \setminus T_Y$.
\end{remark}

\begin{example}
    Let $X$ be a wedge of two triangles and $Y$ a single triangle.
    A simplicial map from $X$ into $Y$ that sends each triangle in $X$ onto the triangle $Y$ is dimension-preserving but not injective.
\end{example}

\begin{corollary}%
    \label{cor:endo-func}
    There is an endofunctor $\tau_1\colon \Simp_{\ord, \dim} \to \Simp_{\ord, \dim}$ defined as follows:
    For $(X, \mathord\preceq_X)$ an ordered simplicial complex, let $T_X$ be the unique order-minimal spanning tree of $X$.
    Then, define $\tau_1(X, \mathord\preceq_X) = (X \boxminus_1 T_X, \mathord\preceq_X)$.

    For a morphism $f\colon (X, \mathord\preceq_X) \to (Y, \mathord\preceq_Y)$ in $\Simp_{\ord, \dim}$, let $T_Y$ be the order-minimal spanning tree of $Y$.
    Define the morphism $\tau_1(f)\colon X \boxminus_1 T_X \to Y \boxminus_1 T_Y$ to be the underlying simplicial map of $f$ formally restricted to $X \boxminus_1 T_X$, interpreted as a morphism in $\Simp_{\ord, \dim}$.

    In particular, this functor is faithful but not full in general.
\end{corollary}

\begin{remark}
    Let $A \subseteq X$ be a simplicial pair and let $i\colon A \to X$ denote the canonical inclusion map.
    Further, we assume that $X$ is endowed with a simplicial order~$\preceq_X$ and we endow $A$ with the restriction of $\preceq_X$ as simplicial order.
    Under this assumption, we see that $i$ is an order-preserving simplicial map $(A, \preceq_X) \to (X, \preceq_X)$.

    Moreover, it follows from Remark~\ref{rem:mono-is-dim-preserv} that $i$ is also dimension-preserving.
    Thus, Corollary~\ref{cor:endo-func} implies that there is an order- and dimension-preserving simplicial map $\tau_1(i)\colon (A \boxminus_1 T_A, \mathord\preceq_X) \to (X \boxminus_1 T_X, \mathord\preceq_X)$, where $T_A$ is the order-minimal spanning tree of~$A$ and $T_X$ is the order-minimal spanning tree of~$X$.

    From this result, it directly follows that $A_1 \setminus T_A$ is a subset of $X_1 \setminus T_X$.
    Hence, Corollary~\ref{cor:endo-func} is a vast generalization of Lemma~\ref{lem:classic-comp-incl} and Theorem~\ref{thm:main}.
\end{remark}

\section{Upper-set decomposables as spans}%
\label{sect:r-spans}

Let $k$ be a principal ideal domain and $Q = (Q, \mathord\leq)$ a ordered set.
Denote by $U_0$ the forgetful functor $\Mod{k} \to \Set$ sending a $k$-module to its underlying set.
The functor $U_0$ is right-adjoint to the $k$-span functor $k\colon \Set \to \Mod{k}$, which assigns to a set $B$ the free $k$-module~$k(B)$ with basis $B$.

Define the functors
\begin{align*}
    U_Q\colon \Mod{kQ} &\longrightarrow \Set^Q \text, &\qquad k_Q\colon \Set^Q &\longrightarrow \Mod{kQ} \text, \\
    M &\longmapsto U_0 \circ M \text, &\qquad B &\longmapsto k \circ B \text.
\end{align*}
It follows from categorical standard arguments that $U_Q$ and $k_Q$ are a pair of adjoint functors.

\begin{lemma}%
    \label{lem:adj}
    There is an adjunction $U_Q \dashv k_Q$.

    \begin{proof}
        Let $\eta = (\eta_{MB}\colon \Hom_{\Mod{k}}(M, k(B)) \to \Hom_{\Set}(U(M), B))_{M,B}$ be the underlying natural isomorphism of the adjunction $U \dashv k$.
        For $N \in \Mod{kQ}$ and $C \in \Set^Q$, define $\eta_{\ast NC}\colon \Hom_{\Mod{kQ}}(N, k_Q(C)) \to \Hom_{\Set^Q}(U_Q(N), C)$ by $\eta_{\ast NC}(\varphi)_q = \eta_{N_q, C_q}(\varphi_q)$ for all $\varphi = (\varphi_q)_{q \in Q}\colon N \to k_Q(C)$ and $q \in Q$.
        It follows from a standard calculation that $\eta_\ast = (\eta_\ast NC)_{N \in \Mod{kQ}, C \in \Set^I}$ is a natural isomorphism from $\Hom_{\Mod{kQ}}(\mathord-, k_Q(\mathord-))$ in $\Hom_{\Set^Q}(U_Q(\mathord-), \mathord-)$.
        Thus, $U_Q$ is left-adjoint to $k_Q$.
    \end{proof}
\end{lemma}

In this sense, the functor $k_Q\colon \Set^Q \to \Mod{kQ}$ can be seen as the persistent generalization of the free $k$-span functor $k\colon \Set \to \Mod{k}$.

\begin{definition}
    Let $Q$ be an ordered set.
    \begin{enumerate}
        \item A subset $U \subseteq Q$ is an \emph{upper set} in $Q$ if
        \[ \{ q \in Q \mid q \geq a \} \subseteq U \]
        for all $a \in U$.
        \item For an upper set~$U$ in $Q$, define the persistent set~$\pt[U]\colon Q \to \Set$ by
        \[ \pt[U](q) = \begin{cases*}
            \{ q \} & if $q \in U$, \\
            \emptyset & otherwise,
        \end{cases*} \]
        with the obvious inclusions.
        The persistent set $\pt[U]$ is called \emph{upper set module} of $U$.
        \item A persistence module~$M\colon Q \to \Mod{k}$ is said to be \emph{upper set decomposable} if there exists a family $(U_i)_{i \in I}$ of upper sets in $Q$ such that
        \[ M \cong \bigoplus_{i \in I} k_Q(\pt[U_i]) \text. \]
    \end{enumerate}
\end{definition}

\begin{proposition}%
    \label{prop:usd}
    Let $B\colon Q \to \Set$ a functor.
    For $q \in Q$, denote by $\iota_q\colon B(q) \to \colim B$ the colimit projection.
    The following assertions are equivalent:
    \begin{enumerate}[label=(\roman*)]
        \item\label{prop:usd:1} For every $q \in Q$, the map $\iota_q$ is injective.
        \item\label{prop:usd:2} The persistence module $\mathop{k_Q} B$ is upper set decomposable.
        \item\label{prop:usd:3} For every $q \in Q$, the colimit projection $\iota'_q\colon \mathop{k_Q} B(q) \to \colim \mathop{k_Q} B$ is a monomorphism.
    \end{enumerate}

    \begin{proof}
        \ref{prop:usd:1} to \ref{prop:usd:2}:
        We define the persistent set $B'\colon Q \to \Set$ as follows:
        For $q \in Q$, let $B'(q) = \iota_q(B(q))$.

        For an ordered pair $q \leq q'$ in $Q$ and every $x \in B'(q)$ we have $x \in B'(q')$ because $\iota_q = \iota_{q'} \circ B(q \leq q')$.
        Therefore, the map $B'(q \leq q')\colon \iota_q(B(q)) \to \iota_{q'}(B(q'))$ can be chosen as the inclusion.

        Moreover, $B'$ is isomorphic to $B$ by construction.
        For every $q \in Q$, the map
        \[ B(q) \to B'(q) \text, \qquad b \mapsto \iota_q(b) \text, \]
        is bijective by the injectivity of $\iota_q$.
        Thus, $B' \cong B$.

        Furthermore, the persistence module~$B'$ can be decomposed into the coproduct
        \[ B' = \coprod_{x \in \colim B} \pt[U_x] \text, \]
        where $U_x = \{ q \in Q \mid x \in B'(q) \}$.
        The set $U_x$ is by construction an upper set.
        Thus, $B'$ is upper set decomposable.

        \ref{prop:usd:2} to \ref{prop:usd:3}:
        Let $(U_i)_{i \in I}$ be a family of upper sets of $Q$ such that $\mathop{k_Q} B = \bigoplus_{i \in I} k[U_i]$.
        Since every upper set module $k[U_i]$ admits injective colimit projections $k[U_i](q) \to \colim k[U_i]$ for any $q \in Q$, it follows that $k_Q(B)$ also admits injective colimit projections as claimed.

        \ref{prop:usd:3} to \ref{prop:usd:1}:
        Suppose that there exists some $q \in Q$ such that $\iota_q$ is not injective.
        Then, there are $b, b' \in B(q)$ such that $b \neq b'$ and $\iota_q(b) = \iota_q(b')$.
        Thus, $b - b' \in k_Q B(q)$ is a non-zero element of the kernel of $\iota'_q$ since $k_Q$ preserves colimits.
        Hence, $\iota'_q$ is not injective and thus not a monomorphism.
    \end{proof}
\end{proposition}

\begin{remark}
    A similar result was established in \cite[Proposition~3.1]{ChachólskiScolamieroVaccarino2017}.
    A functor $B\colon Q \to \Mod{k}$ is called \emph{multifiltration} if the structure maps of $B$ are monomorphisms.
    For the ordered set~$Q = \mathbb N^d$ (or more generally for $Q$ a lattice), the conditions of Proposition~\ref{prop:usd} are equivalent to the condition that the functor $B$ is a multifiltration, cf. \cite[Proposition~3.3]{ChachólskiScolamieroVaccarino2017}.
    Moreover, any non-trivial multifiltration $B$ is uniquely decomposable into indecomposables and the indecomposable multifiltrations over $\mathbb N^d$ are precisely the $k$-linear generators of upper set modules \cite[Proposition~3.1, Corollary~3.2]{ChachólskiScolamieroVaccarino2017}.
    In this sense, Proposition~\ref{prop:usd} is a generalization of \cite[Proposition~3.1]{ChachólskiScolamieroVaccarino2017}.
\end{remark}

\begin{example}
    Let $U = \{ (a, b) \in \mathbb Z^2 \mid a + b \geq 0 \}$.
    This is clearly an upper set of $\mathbb Z^2$.
    However, the upper set module $k[U]$ is not finitely generated, since we need a generator in each $(a, -a)$ for $a \in \mathbb Z$ to cover $U$ by principal upper sets of $\mathbb Z^2$.
\end{example}

\section{Upper set generators of 1-dimensional persistent homology}
\label{sect:gen-h1}

Let $X = (X^q)_{q \in Q}$ be a $Q$-filtration of simplicial complexes, and let $T = (T^q)_{q \in Q}$ be a subfiltration of spanning trees for $X$.
By Proposition~\ref{prop:pers-filt}, the quotient persistence module $C_1(X, T) = C_1(X) / C_1(T)$ is isomorphic to $Z_1(X)$.

Because $C_1(X, T)$ has injective colimit projections and a persistent basis, Proposition~\ref{prop:usd} implies that $C_1(X, T)$ is upper set decomposable.
Therefore, there is a decomposition $Z_1(X) \cong \bigoplus_{i \in I} k[U_i]$ for some family $(U_i)_{i \in I}$ of upper sets in $Q$.

However, the existence of this decomposition depends on the requirement that $T$ is a subfiltration of spanning trees.
Subfiltrations of spanning trees do not always exist, as demonstrated in Example~\ref{ex:no-subfil}, if $Q$ is an arbitrarily ordered set.
In fact, even for a $3 \times 3$-sublattice of $\mathbb N^2$, the existence of filtrations of spanning trees fails.

As an alternative for the general case, we construct a covering of $Z_1(X)$ by an upper set decomposable module, with a cofiltration of spanning trees of $X$ as the main input.
Explicitly, we choose the cofiltration of order-minimal spanning trees, which admits the properties in Theorem~\ref{thm:main}.

We assume that the ordered set~$Q$ and the simplicial complex $\bigcup_{q \in Q} X^q$ are finite and endowed with a simplicial order~$\preceq$.
Then, for every $q \in Q$ we interpret the subcomplex $X^q$ as the ordered simplicial complex $(X^q, \mathord\preceq)$, which naturally admits an injective order-preserving simplicial map into $(X, \mathord\preceq)$.
We choose the uniquely determined order-minimal spanning tree $T^q$ of $(X^q, \mathord\preceq)$, i.e.\ $T^q = \tau_1(X^q, \mathord\preceq)$.

It follows from Theorem~\ref{thm:main} that $X^q \boxminus_1 T^q$ is a subcomplex of $X^{q'} \boxminus_1 T^{q'}$ for every ordered pair $q \leq q'$ in $Q$.
Thus, $T = (T^q)_{q \in Q} = (\tau_1(X^q))_{q \in Q}$ is a cofiltration of spanning trees of $X$.

\begin{definition}
    Let $\sigma$ be an edge in the total simplicial complex $\bigcup_{q \in Q} X^q \boxminus_1 T^q$.
    The \emph{representative persistent set} of $\sigma$ is
    \[ \mathcal B_\sigma^q = \{ j^{-1}_{q'}(\sigma) \mid q' \in Q \text,\ q' \leq q \text,\ \sigma \in X^{q'} \boxminus_1 T^{q'} \} \text, \]
    where $j_{q'}\colon Z_1(X^{q'}) \to C_1(X^{q'}, T^{q'})$ is the isomorphism defined by $j_{q'}(z) = z + C_1(T^{q'})$.
\end{definition}

\begin{remark}
    Since $(T^q)_{q \in Q}$ is a cofiltration of $X$, it follows that the structure maps of the persistent set $\mathcal B_\sigma$ are given by inclusions.
\end{remark}

\begin{definition}
    The \emph{associated upper set precover} of $X = (X^q)_{q \in Q}$ with respect to the cofiltration of spanning trees~$T$ is defined by
    \[ \ug_T(X) = \bigoplus_{\sigma} k_Q \mathcal B_\sigma \text, \]
    where $\sigma$ goes over all edges in $\bigcup_{q \in Q} X^q \boxminus_1 T^q$.
\end{definition}

\begin{remark}
    It follows from the adjunction $U_Q \vdash k_Q$, Lemma~\ref{lem:adj}, that $\ug_T(X) = k_Q(\coprod_\sigma B_\sigma)$, where $\sigma$ goes over all edges in $\bigcup_{q \in Q} X^q \boxminus_1 T^q$.
\end{remark}

\begin{theorem}%
    \label{thm:gen-h1}
    The associated upper set precover $\ug_T(X)$ is upper set decomposable and admits an epimorphism $\ug_T(X) \to Z_1(X)$.

    \begin{proof}
        The persistent set $\mathcal B = (\mathcal B^q)_{q \in Q} = \coprod_{\sigma} \mathcal B_\sigma$ is given by
        \[ \mathcal B^q = \{ (\sigma, j^{-1}_{q'}(\sigma)) \mid q' \in Q\text,\ q' \leq q\text,\ \sigma \in (X^q \boxminus_1 T^q) \cap (X^{q'} \boxminus_1 T^{q'}) \} \text. \]
        We define the persistent map $\varphi = (\varphi^q)_{q \in Q} \colon \mathcal B \to \mathop{U_Q} Z_1(X)$ as follows:
        For $(\sigma, z) \in \mathcal B^q$, $z$~is a cycle of $X^{q'}$ for some $q' \in Q$ satisfying $q' \leq q$.
        Then, $z$~is a cycle in $X^q$ because $X^{q'} \subseteq X^{q}$, and we define $\varphi^q(\sigma, z) = z \in Z_1(X^q)$.
        This induces a uniquely determined homomorphism $\tilde\varphi^q\colon \ug_T(X) = k_Q(\mathcal B) \to Z_1(X^q)$.

        Furthermore, for every $q \in Q$ we have a basis
        \[ B_q = \{ j_q^{-1}(\sigma + C_1(T^q)) \mid \sigma \in (X^q \boxminus_1 T^q)_1 \} \]
        of $Z_1(X^q)$, which is contained in $\mathcal B^q$.
        Thus, the image of $\tilde\varphi^q$ contains $B_q$, and $\tilde\varphi^q$ is pointwise surjective.
        Hence, $\tilde\varphi$ is an epimorphism in the category of persistence modules.
    \end{proof}
\end{theorem}

\section{Higher-dimensional spanning complexes}%
\label{sect:n-span}

In the previous sections we have discussed spanning trees of simplicial complexes, which allow us to determine the $1$-dimensional cycle modules of simplicial complexes by inspecting $1$-dimensional relative chain modules.
An open question is whether there is a similar, easily computable type of subcomplex such that for a dimension~$n \in \mathbb N$ we can determine $n$-dimensional cycle modules in a similar manner by inspecting $n$-dimensional relative chain modules.

As a first step towards an answer of this question, we consider in this section a generalization of spanning trees to higher dimensions by merely generalizing Definition~\ref{def:st}.
\begin{definition}%
    \label{def:n-span}
    Let $X$ be a simplicial complex and $n > 0$.
    A subcomplex~$A$ of~$X$ is an \emph{$n$-spanning complex} of~$X$ if the following conditions hold:
    \begin{enumerate}[label=(\roman*)]
        \item\label{def:n-span:1} $\skel_{n-1}(A) = \skel_{n-1}(X)$.
        \item\label{def:n-span:2} The boundary homomorphism $\partial_n\colon C_n(A) \to C_{n-1}(A)$ is injective.
        \item\label{def:n-span:3} $B_{n-1}(A) = B_{n-1}(X)$.
    \end{enumerate}
\end{definition}

This definition directly generalizes the notion of a spanning tree defined in Definition~\ref{def:st} by replacing the dimensions by an fixed but arbitrary $n \in \mathbb N$.

\begin{remark}
    For $n = 1$, every spanning tree is an $n$-spanning complex.
\end{remark}

\begin{proposition}%
    \label{prop:n-span-main}
    Let $A$ be an $n$-spanning complex of $X$.
    The homomorphism
    \begin{align*}
        j\colon Z_n(A) &\longrightarrow C_n(X, A) \text, \\
        z &\longmapsto z + C_n(A) \text,
    \end{align*}
    is an isomorphism of $k$-modules.

    \begin{proof}
        The proof is analogous to Proposition~\ref{prop:st-main}.
        The injectivity of $j$ follows from $Z_n(A) = 0$ and the surjectivity of $j$ follows from $B_{n-1}(X) = B_{n-1}(A)$ together with $Z_n(A) = 0$.
    \end{proof}
\end{proposition}

\begin{theorem}[Existence of $n$-spanning complexes]%
    \label{thm:exist-n-span}
    For any simplicial complex~$X$ and any $n \in \mathbb N$ there exists a subcomplex $A$ such that $Z_n(A) = 0$ and $B_{n-1}(A) = B_{n-1}(X)$, i.e.\ $X$ admits an $n$-spanning complex.

    \begin{proof}
        Let $\mathcal A_X = \{ A \subseteq X \mid \text{$A$ is a subcomplex s.t.\ $Z_n(A) = 0$} \}$.
        For every subset $\mathcal A' \subseteq \mathcal A_X$ that is totally ordered by set inclusion, we consider the subcomplex $\bigcup \mathcal A = \bigcup_{A \in \mathcal A'} A$ of $X$.
        We have $Z_n(\bigcup \mathcal A') = 0$:
        Let $z$ be an $n$-chain in $\bigcup \mathcal A'$ such that $\partial_n z = 0$.
        Then, there are non-zero $\lambda_1, \dotsc, \lambda_m \in k$ and distinct $n$-simplices $\sigma_1, \dotsc, \sigma_m$ in $\bigcup \mathcal A'$ such that $z = \lambda_1 \sigma_1 + \dotsb + \lambda_m \sigma_m$.
        Then, we can choose $A_1, \dotsc, A_m \in \mathcal A'$ such that $\sigma_i \in A_i$ for $i \in \{ 1, \dotsc, m \}$.
        Since $\mathcal A'$ is totally ordered by set inclusion there exists an $i_0 \in \{ 1, \dotsc, m \}$ such that $A_i \subseteq A_{i_0}$ for all $i \in \{ 1, \dotsc, m \}$.
        Thus, $z$ is actually an $n$-chain in $B_{i_0} \in \mathcal A'$.
        Since $B_{i_0}$ satisfies $Z_1(B_{i_0}) = 0$ by assumption, $z = 0$.

        Therefore, every totally ordered subset of $\mathcal A_X$ admits an upper bound in $\mathcal A_X$.
        It follows from Zorn's lemma that there exists $A \in \mathcal A_X$ that is maximal w.r.t.\ set inclusion.

        It remains to prove that such maximal $A$ in $\mathcal A_X$ satisfies $B_{n-1}(A) = B_{n-1}(X)$, i.e.\ $B_{n-1}(X) \subseteq B_{n-1}(A)$.
        In particular, it suffices to prove that for every $n$-simplex $\sigma$ in $X$ there exists an $n$-chain $c$ in $A$ such that $\partial_n c = \partial_n \sigma$.
        Let $\sigma$ be an $n$-simplex in $X$.
        If $\sigma \in A$, then trivially $\partial_n \sigma \in B_{n-1}(A)$.
        Thus, we consider only the case that $\sigma$ is not in $A$.
        Then, we can consider the subcomplex $A' = A \cup \{ \sigma \}$ of $X$.
        Since $A$ is maximal by construction in $\mathcal A_X$ and $A \subsetneq A'$, it follows that $A'$ cannot be contained in $\mathcal A_X$.
        That is, $Z_n(A')$ cannot be trivial.

        Therefore, there exists an $n$-cycle $z = \mu \sigma + c$ in $A'$ with $\mu \in \{ 1, -1 \}$ and $c$ an $n$-chain in $A$.
        We can assume without loss of generality that $\mu = 1$.
        But then we have $\partial_n \sigma = \partial_n c \in B_{n-1}(A)$.
        Hence, $B_{n-1}(X) \subseteq B_{n-1}(A)$.
    \end{proof}
\end{theorem}

{\emergencystretch=1em
\printbibliography}
\end{document}